\documentclass[11pt]{article}

\usepackage{a4wide}
\usepackage{amssymb}
\usepackage{amsthm}
\usepackage{enumerate}
\usepackage{amsmath}

\newtheorem{thm}{Theorem}
\newtheorem*{defi}{Definition}
\newtheorem{cor}[thm]{Corollary}

\newtheorem{prop}[thm]{Proposition}

\newcommand{\N}{\mathbb N}
\renewcommand{\L}{\mathbb L}
\newcommand{\Z}{\mathbb Z}
\newcommand{\R}{\mathbb R}
\newcommand{\dif}{{\rm d}}
\newcommand{\ind}{{\mathchoice{\mathrm {1\mskip-4.1mu l}} {\mathrm{ 1\mskip-4.1mu l}} {\mathrm {1\mskip-4.6mu l}} {\mathrm {1 \mskip-5.2mu l}}}}

\newcommand{\F}{{\mathcal F}}
\newcommand{\A}{{\mathcal A}}
\newcommand{\B}{{\mathcal B}}

\newcommand{\esp}{\mathbb E}

\newcommand{\cqfd}{\hfill $\Box$}

\begin{document}

 \title{\Huge{A fourth moment inequality for\\ functionals of stationary processes}}
 \author{ Olivier Durieu\\
\\
Laboratoire de Math\'ematique Rapha\"el Salem,\\ UMR 6085 CNRS-Universit\'e de Rouen,\\
e-mail: olivier.durieu@etu.univ-rouen.fr}

\maketitle

\begin{abstract}
In this paper, a fourth moment bound for partial sums of functional of strongly ergodic Markov chain is established.
This type of inequality plays an important role in the study of empirical process invariance principle. This one is 
specially adapted to the technique of Dehling, Durieu and Voln\'y (2008).
The same moment bound can be proved for dynamical system whose transfer operator has some spectral properties.
Examples of applications are given.

\medskip

\noindent Keywords: Stationary processes; Moment inequalities; Strongly ergodic Markov chains; Dynamical system; Empirical distribution; Invariance principle.

\medskip

\noindent AMS classification: 60G10; 60J10; 60F17; 28D05; 62G20.
\end{abstract}

\section{Introduction}\label{s1}

Fourth moment bounds for partial sums of stationary processes are a key tool in the study of functional limit theorems.
In particular, it is the case in the investigation of empirical process invariance principle. It will be our main application.
Let $(X_i)_{i\ge 0}$ be a stationary process in $\R$ and $F(t)=P(X_0\leq t)$. The empirical process associated to $(X_i)_{i\ge 0}$ is defined by
\begin{equation}
Y_n(t)=\frac{1}{\sqrt{n}}\sum_{i=0}^{n-1}\{\ind_{(-\infty,t]}(X_i)-F(t)\},\quad t\in\R.
\end{equation}
One says that $(Y_n)_{n\ge 0}$ satisfies an invariance principle if it converges in distribution to a mean-zero Gaussian process.
In general, proofs of invariance principle consist of two parts: multivariate central limit theorem and tightness.
Since the work of Donsker \cite{Don52} for i.i.d. sequences and the work of Billingsley \cite{Bil68} for some weakly dependent
processes, the chaining technique seems to be a suitable way to get the tightness of the process $(Y_n)_{n\ge0}$.
Fourth moment inequalities are a central point in this technique.
In many cases (like i.i.d. case), if an inequality of the type
\begin{equation}\label{mb}
\esp\left[\sum_{i=0}^{n-1}\{\ind_{(s,t]}(X_i) -(F(t)-F(s))\}\right]^4\leq C(n(t-s)+n^2(t-s)^2),
\end{equation}
for all $s<t$ is established then the tightness of the process follows.
The difficulty is to deal with the sequence of indicator variables $(\ind_{(s,t]}(X_i))_{i\ge 0}$. When the process
$(X_i)_{i\ge 0}$ has some mixing properties, like strong, uniform or beta mixing, $(\ind_{(s,t]}(X_i))_{i\ge 0}$
inherits comparable properties and fourth moment bound can be established assuming some regularity conditions on $X_0$'s distribution function (for an overview of this theory, see
Dehling and Philipp \cite{DehPhi02} and references therein). In 2004, Collet, Martinez and Schmitt \cite{ColMarSch04}
proved a fourth moment inequality for expanding maps of the interval (see also Dedecker and Prieur \cite{DedPri07}).
They use spectral properties of the transfer
operator associated on the space of bounded variation functions and the fact that indicators belong to this space. 
In other cases, some type of Markov chains and dynamical system (see Hennion and Herv\'e \cite{HenHer01}, Gou\"ezel \cite{Gou08}),
one can have good properties
of the Markov operator or the transfer operator on other spaces of functions which do not contain the indicator functions.
In Dehling, Durieu and Voln\'y \cite{DehDurVol08}, a new technique for proving empirical process invariance principle is developed using approximations of indicators by regular functions. The fourth moment bound of Corollary \ref{cor} is well adapted to this situation. The main point is that the Banach norm only appears through its logarithm. In Section \ref{s2} and Section \ref{s3} the fourth moment inequality is stated and proved for strongly ergodic Markov chain. In Section \ref{s4}, we state the same moment bound
for a class of dynamical systems. Section \ref{s5} is devoted to examples of applications.

\section{Fourth moment inequality for Markov chains}\label{s2}
Let $\left(X_n\right)_{n\ge 0}$ be a homogeneous Markov chain with a stationary measure $\nu$.
Denote by $P$ the associated Markov operator and $E$ the state space.
Consider a Banach space $\left(\B,\|. \|\right)$ of $\nu$-measurable functions from $E$ to $\R$, which contains
the function $\bf{1}=\ind_E$.
We will assume that the chain is $\B$-geometrically ergodic.

\begin{defi}
The Markov chain $\left(X_n\right)_{n\ge 0}$ is $\B$-geometrically ergodic or strongly ergodic (with respect to $\B$)
if
\begin{enumerate}
\item[a)] there exist $\kappa>0$ and $0<\theta<1$ such that for all $f\in\B$, 
\begin{equation}\label{exp}
\|P^nf-\Pi f\|\leq \kappa\theta^n\|f\|
\end{equation}
where $\Pi f=\esp_\nu\left(f\right)\bf{1}$.

\item[b)] there exists $p\geq 1$ such that $\left(\B,\|. \|\right)$ is continuously included in $\left(\L^p\left(\nu\right), \|.\|_p\right)$.\\
i.e. $\exists C>0$ such that $\forall f\in\B,$ 
\begin{equation}\label{cont}
\|f\|_p\leq C\|f\|
\end{equation}
where $\|f\|_p=\left(\int |f|^p\dif\nu\right)^\frac{1}{p}$.
\end{enumerate}
\end{defi}

\medskip

Further, we assume that there exists a constant $M>0$ such that $\forall f\in\B$ and $\forall n\in\N$,
\begin{equation}\label{algebra}
\|fP^nf\|\leq M\|f\|\|P^nf\|.
\end{equation}
In particular, this is the case if $(\B,\|.\|)$ is a Banach algebra. In the sequel, with no loss of generality, we assume $M=1$.

Strong ergodicity covers a large class of examples (discussed in Section 5). It corresponds to the fact that the Markov
transition operator acting on $\B$ has $1$ as simple eigenvalue and the rest of the spectrum is included in a closed ball of radius strictly smaller than $1$.

\medskip

For a function $\varphi:E\longrightarrow\R$, we consider the partial sum 
$$
S_n\left(\varphi\right)=\sum_{i=1}^n\varphi\left(X_i\right).
$$

The aim is to get a fourth moment inequality for this partial sum when the function $\varphi$ belongs to the space $\B$.
Our main results are the following.

\begin{thm}\label{thm}
If $\left(X_n\right)_{n\ge 0}$ is a $\B$-geometrically ergodic Markov chain with stationary measure $\nu$, then
for all $\varphi\in\B$ such that $\esp_\nu\left(\varphi\right)=0$ and $\varphi\in\L^4(\nu)\cap\L^{3q}\left(\nu\right)$,
\begin{eqnarray*}
&&\esp_\nu\left[S_n\left(\varphi\right)^4\right]\\
&&\leq K\left[n\|\varphi\left(X_0\right)^4\|_1\log^3(\|\varphi\|+1)\right.\\
&&\quad\quad +n\left(\|\varphi\left(X_0\right)^3\|_q+\|\varphi\left(X_0\right)^2\|_q+\|\varphi\left(X_0\right)\|_q+
\|\varphi\left(X_0\right)\|_q^2\right)\log^2(\|\varphi\|+1)\\
&&\left.\quad\quad+n^2\left(\|\varphi\left(X_0\right)^2\|_1\log(\|\varphi\|+1)+\|\varphi\left(X_0\right)\|_q\right)^2\right]
\end{eqnarray*}
where $\frac{1}{p}+\frac{1}{q}=1$ and $K$ is a constant.
\end{thm}

\medskip

As consequence, the following corollary give a simplest inequality is the case where the function $\varphi$ is bounded.

\begin{cor}\label{cor}
If $\left(X_n\right)_{n\ge 0}$ is a $\B$-geometrically ergodic Markov chain with stationary measure $\nu$, then
for all $\varphi\in\B$ such that $\esp_\nu\left(\varphi\right)=0$, $\varphi$ is bounded and
${\rm m}_\varphi=\max\{1,\sup_x|\varphi(x)|\}$,
$$
\esp_\nu\left[S_n\left(\varphi\right)^4\right]\leq
K{\rm m}_\varphi^3\left[n\|\varphi\left(X_0\right)\|_q\log^3(\|\varphi\|+1)
+n^2\|\varphi\left(X_0\right)\|_q^2\log^2(\|\varphi\|+1)\right].
$$
\end{cor}

\bigskip

Assume that one can prove a multivariate central limit theorem for functions in $\B$, then by the technique of Dehling, Durieu and Voln\'y \cite{DehDurVol08}, if the space $\B$ contains a class of functions approximating the indicators, an empirical process invariance principle  follows.

To be complete, we state the following result, which is a corollary of Gordin's Theorem \cite{Gor69} and from which we can deduce a multivariate central limit theorem.
\begin{prop}
If $(X_n)_{n\ge 0}$ is ergodic and $\B$-geometrically ergodic with $p\ge 2$ in (\ref{cont}), then
for all $\varphi\in\B$, $\frac{S_n(\varphi-\esp_\nu\varphi)}{\sqrt{n}}$ converges in distribution to a centred normal law.
\end{prop}
See Hennion and Herv\'e \cite{HenHer01}, Thm.IX.2, for sufficient conditions on $\B$ to have: strong ergodicity implies ergodicity.

\section{Proof of Theorem \ref{thm}}\label{s3}
Let us consider the assumptions of Theorem \ref{thm} hold.
In the sequel, all the expectations are considered with respect to the measure $\nu$
 and $\F_i$ denotes the $\sigma$-algebra generated by $X_i$.

Let $\varphi\in\B$ with $\esp\left(\varphi\right)=0$.
As the Markov chain is a stationary process, we have the following bound:
$$
\esp\left[S_n\left(\varphi\right)^4\right]\leq 4!n\sum_{\begin{array}{c}_{i,j,k\ge 0:}\\_{i+j+k\leq n}\end{array}}
\esp\left(\varphi\left(X_0\right)\varphi\left(X_i\right)\varphi\left(X_{i+j}\right)\varphi\left(X_{i+j+k}\right)\right)
$$
which can be decomposed in three sums.
\begin{eqnarray}
\esp\left[S_n\left(\varphi\right)^4\right]&\leq
&4!n\left[\,\sum_{i=1}^n\sum_{j,k\leq i}\esp\left(\varphi\left(X_0\right)\varphi\left(X_i\right)
\varphi\left(X_{i+j}\right)\varphi\left(X_{i+j+k}\right)\right)\right.\label{i}\\
&&\qquad+\sum_{j=1}^n\sum_{i,k\leq j}\esp\left(\varphi\left(X_0\right)\varphi\left(X_i\right)
\varphi\left(X_{i+j}\right)\varphi\left(X_{i+j+k}\right)\right)\label{j}\\
&&\qquad+\left.\sum_{k=1}^n\sum_{i,j\leq k}\esp\left(\varphi\left(X_0\right)\varphi\left(X_i\right)
\varphi\left(X_{i+j}\right)\varphi\left(X_{i+j+k}\right)\right)\right].\label{k}
\end{eqnarray}

So to study the terms $\esp\left(\varphi\left(X_0\right)\varphi\left(X_i\right)\varphi\left(X_{i+j}\right)\varphi\left(X_{i+j+k}\right)\right)$,
we will consider three cases according to the greatest integer between $i, j$ and $k$.

\medskip

Before, we can see that in all cases, by H\"older inequality, we have
\begin{equation}
\left|\esp\left(\varphi\left(X_0\right)\varphi\left(X_i\right)\varphi\left(X_{i+j}\right)\varphi\left(X_{i+j+k}\right)\right)\right|
\leq \|\varphi\left(X_0\right)\|_{4}^4\leq \|\varphi\left(X_0\right)^4\|_{1}.\label{cas1.1}
\end{equation}

Further, let $n_0$ be a positive integer such that 
$$
\frac{\log(\|\varphi\|+1)}{-\log\theta}< n_0\leq \frac{\log(\|\varphi\|+1)}{-\log\theta}+1.
$$
Note that $\theta^{n_0}\|\varphi\|\leq 1$. 

\bigskip

\textsl{First case}: $i, j \leq k$.

\medskip

Here, we  use the properties of the Markov operator $P$ on the space $\B$ to get another majoration.
Applying successively H\"older inequality, properties (\ref{cont}) and (\ref{exp}), we get
\begin{eqnarray}
&&|\esp\left(\varphi\left(X_0\right)\varphi\left(X_i\right)\varphi\left(X_{i+j}\right)\varphi\left(X_{i+j+k}\right)\right)|\nonumber\\
&&=|\esp\left(\varphi\left(X_0\right)\varphi\left(X_i\right)\varphi\left(X_{i+j}\right)\left(\esp\left(\varphi\left(X_{i+j+k}\right)|\F_{i+j}\right)-\esp\left(\varphi\left(X_{i+j+k}\right)\right)\right)\right)|\nonumber\\
&&\leq\|\varphi\left(X_0\right)\varphi\left(X_i\right)\varphi\left(X_{i+j}\right)\|_q\|P^k\varphi\left(X_0\right)-\Pi\varphi\left(X_0\right)\|_p\nonumber\\
&&\leq \|\varphi\left(X_0\right)\|_{3q}^3 C\|P^k\varphi-\Pi\varphi\|\nonumber\\
&&\leq \|\varphi\left(X_0\right)^3\|_q C\kappa\theta^k\|\varphi\|.\label{cas1.2}
\end{eqnarray}

\medskip

Now, for sum $\left(\ref{k}\right)$, using $\left(\ref{cas1.1}\right)$ for the $n_0-1$ first terms and $\left(\ref{cas1.2}\right)$ for the others, we obtain
\begin{eqnarray*}
&&\sum_{k=1}^n\sum_{i,j\leq k}\esp\left(\varphi\left(X_0\right)\varphi\left(X_i\right)\varphi\left(X_{i+j}\right)\varphi\left(X_{i+j+k}\right)\right)\\
&&\leq \sum_{k=1}^{n_0-1}k^2\|\varphi\left(X_0\right)^4\|_1+\sum_{k=n_0}^nk^2C\kappa\|\varphi\left(X_0\right)^3\|_q\theta^k\|\varphi\|\\
&&\leq (n_0-1)^3\|\varphi\left(X_0\right)^4\|_1+C\kappa\|\varphi\left(X_0\right)^3\|_q\sum_{k=n_0}^nk^2\theta^{k-n_0}.
\end{eqnarray*}

There exists a constant $C_1$ which only depends on $\theta$, such that
\begin{eqnarray*}
\sum_{k=n_0}^nk^2\theta^{k-n_0}&\leq&\sum_{k\geq 2}\left(k+n_0-2\right)^2\theta^{k-2}\\
&\leq&\sum_{k\geq 2}k^2\theta^{k-2}+2(n_0-2)\sum_{k\geq 2}k\theta^{k-2}+(n_0-2)^2\sum_{k\geq 2}\theta^{k-2}k\\
&\leq& C_1(n_0-1)^2,
\end{eqnarray*}
because the three series converge.

\medskip

Thus, writing $C_2=-\frac{1}{\log\theta}$, we get $n_0-1\leq C_2\log(\|\varphi\|+1)$ and 
\begin{eqnarray}
&&\sum_{k=1}^n\sum_{i,j\leq k}\esp\left(\varphi\left(X_0\right)\varphi\left(X_i\right)\varphi\left(X_{i+j}\right)\varphi\left(X_{i+j+k}\right)\right)\nonumber\\
&&\leq C_2^3\|\varphi\left(X_0\right)^4\|_1\log^3(\|\varphi\|+1)
+C_3\|\varphi\left(X_0\right)^3\|_q\log^2(\|\varphi\|+1),\label{cas1}
\end{eqnarray}
where $C_3=C\kappa C_1C_2^2$.
\bigskip

\textsl{Second case}: $i, k \leq j$.

\medskip

We can decompose the expectation as follows,
\begin{eqnarray}
&&|\esp\left(\varphi\left(X_0\right)\varphi\left(X_i\right)\varphi\left(X_{i+j}\right)\varphi\left(X_{i+j+k}\right)\right)|\nonumber\\
&&\leq|\esp\left(\varphi\left(X_0\right)\varphi\left(X_i\right)\left(\esp\left(\varphi\left(X_{i+j}\right)\esp\left(\varphi\left(X_{i+j+k}\right)|\F_{i+j}\right)|\F_{i}\right)-\esp\left(\varphi\left(X_0\right)\varphi\left(X_{k}\right)\right)\right)\right)|\nonumber\\
&&\quad+|\esp\left(\varphi\left(X_0\right)\varphi\left(X_i\right)\right)\esp\left(\varphi\left(X_0\right)\varphi\left(X_k\right)\right)|.\label{cas2.0}
\end{eqnarray}
In the right hand side, let us call $I_{i,j,k}$ the first term and $II_{i,k}$ the second one.

\medskip

Since $\esp\left(\varphi P^k\varphi \left(X_0\right)\right)=\esp\left(\varphi\left(X_0\right)\esp\left(\varphi \left(X_k\right)|\F_0\right)\right)=\esp\left(\varphi\left(X_0\right)\varphi\left(X_k\right)\right)$, we have
\begin{eqnarray}
I_{i,j,k}&\leq &\|\varphi\left(X_0\right)\varphi\left(X_i\right)\|_q\|P^j\left(\varphi P^k\varphi\right)\left(X_0\right)-\Pi\left(\varphi P^k\varphi\right)\left(X_0\right)\|_p\nonumber\\
&\leq & \|\varphi\left(X_0\right)\|_{2q}^2 C\|P^j\left(\varphi P^k\varphi\right)-\Pi\left(\varphi P^k\varphi\right)\|\nonumber\\
&\leq & C\|\varphi\left(X_0\right)^2\|_q\kappa\theta^j\|\varphi P^k\varphi\|\nonumber
\end{eqnarray}
and by assumption (\ref{algebra}),
$$
\|\varphi P^k\varphi\|\leq \|\varphi\| \|P^k\varphi\| \leq \kappa\theta^k\|\varphi\|^2.
$$
Therefore,
\begin{eqnarray*}
I_{i,j,k}&\leq& C\kappa^2\|\varphi\left(X_0\right)^2\|_q\theta^{j+k}\|\varphi\|^2 \\
&\leq& C\kappa^2\|\varphi\left(X_0\right)^2\|_q\theta^{j}\|\varphi\|^2.\label{cas2.1}
\end{eqnarray*}

\medskip

Now, thanks to the decomposition $\left(\ref{cas2.0}\right)$ (using also inequality $\left(\ref{cas1.1}\right)$), for $n$ big enough,
\begin{eqnarray}
&&\sum_{j=1}^n\sum_{i,k\leq j}\esp\left(\varphi\left(X_0\right)\varphi\left(X_i\right)\varphi\left(X_{i+j}\right)\varphi\left(X_{i+j+k}\right)\right)\nonumber\\
&&\leq \sum_{j=1}^{2n_0-2}j^2\|\varphi\left(X_0\right)^4\|_1+\sum_{j=2n_0-1}^n\sum_{i,k\leq j}\left(I_{i,j,k}+II_{i,k}\right)\nonumber\\
&&\leq 8\left(n_0-1\right)^3\|\varphi\left(X_0\right)^4\|_1+\sum_{j=2n_0-1}^nj^2I_{i,j,k}+\sum_{j=1}^n\sum_{i,k\leq j}II_{i,k},\nonumber
\end{eqnarray}
where $n_0$ has been defined previously.

\medskip

Inequality $\left(\ref{cas2.1}\right)$ and $\theta^{n_0}\|\varphi\|\leq 1$ imply
\begin{eqnarray}
\sum_{j=2n_0-1}^nj^2I_{i,j,k}&\leq &\sum_{j=2n_0-1}^nj^2C\kappa^2\|\varphi\left(X_0\right)^2\|_q\theta^{j}\|\varphi\|^2\nonumber\\
&\leq& C\kappa^2\|\varphi\left(X_0\right)^2\|_q\sum_{j=2n_0-1}^nj^2\theta^{j-2n_0}.\nonumber
\end{eqnarray}

As before, there exists a constant $C_4$ depending on $\theta$ such that, 
$$ \sum_{j=2n_0-1}^nj^2\theta^{j-2n_0} \leq C_4(n_0-1)^2. $$

So,
$$
\sum_{j=2n_0-1}^nj^2I_{i,j,k}\leq C\kappa^2\|\varphi\left(X_0\right)^2\|_q C_4(n_0-1)^2.
$$

\medskip

For the third term, we have
$$
\sum_{j=1}^n\sum_{i,k\leq j}II_{i,k}\leq n\left(\sum_{i=1}^n|\esp\left(\varphi\left(X_0\right)\varphi\left(X_i\right)\right)|\right)\left(\sum_{k=1}^n|\esp\left(\varphi\left(X_0\right)\varphi\left(X_k\right)\right)|\right).
$$
We can see that
\begin{eqnarray}
|\esp\left(\varphi\left(X_0\right)\varphi\left(X_i\right)\right)|&\leq& \|\varphi\left(X_0\right)\|_q\|P^i\varphi\left(X_0\right)-\Pi\varphi\left(X_0\right)\|_p\nonumber\\
&\leq& C\|\varphi\left(X_0\right)\|_q\|P^i\varphi-\Pi\varphi\|\nonumber\\
&\leq& C\kappa\|\varphi\left(X_0\right)\|_q\theta^i\|\varphi\|\label{cas2.2}
\end{eqnarray}
and, in the same way,
\begin{equation}
|\esp\left(\varphi\left(X_0\right)\varphi\left(X_k\right)\right)|\leq C\kappa\|\varphi\left(X_0\right)\|_q\theta^k\|\varphi\|.\label{cas2.3}
\end{equation}
Alternatively, by H\"older inequality, 
\begin{equation}
|\esp\left(\varphi\left(X_0\right)\varphi\left(X_i\right)\right)|\leq \|\varphi\left(X_0\right)^2\|_1\quad{\rm and}\quad |\esp\left(\varphi\left(X_0\right)\varphi\left(X_k\right)\right)|\leq \|\varphi\left(X_0\right)^2\|_1.\label{cas2.4}
\end{equation}

Thus, by $\left(\ref{cas2.2}\right)$, $\left(\ref{cas2.3}\right)$ and $\left(\ref{cas2.4}\right)$,
\begin{eqnarray}
\sum_{i=1}^n|\esp\left(\varphi\left(X_0\right)\varphi\left(X_i\right)\right)|&\leq& \sum_{i=1}^{n_0-1}\|\varphi\left(X_0\right)^2\|_1+\sum_{i=n_0}^nC\kappa\|\varphi\left(X_0\right)\|_q\theta^i\|\varphi\|\nonumber\\
&\leq& (n_0-1)\|\varphi\left(X_0\right)^2\|_1+C\kappa\|\varphi\left(X_0\right)\|_q\sum_{i=n_0}^n\theta^{i-n_0}\nonumber\\
&\leq& (n_0-1)\|\varphi\left(X_0\right)^2\|_1+C_5\|\varphi\left(X_0\right)\|_q,\nonumber
\end{eqnarray}
where $C_5=C\kappa\sum_{i\geq 0}\theta^i<\infty$ and
\begin{equation}
\sum_{k=1}^n|\esp\left(\varphi\left(X_0\right)\varphi\left(X_k\right)\right)|\leq (n_0-1)\|\varphi\left(X_0\right)^2\|_1+C_5\|\varphi\left(X_0\right)\|_q.\nonumber
\end{equation}

\medskip

Finally,
\begin{eqnarray}
&&\sum_{j=1}^n\sum_{i,k\leq j}\esp\left(\varphi\left(X_0\right)\varphi\left(X_i\right)\varphi\left(X_{i+j}\right)\varphi\left(X_{i+j+k}\right)\right)\nonumber\\
&&\leq 8C_2^3\|\varphi\left(X_0\right)^4\|_1\log^3(\|\varphi\|+1)+ C_6\|\varphi\left(X_0\right)^2\|_q\log^2(\|\varphi\|+1)\label{cas2}\\
&&\quad + n\left(C_2\|\varphi\left(X_0\right)^2\|_1\log(\|\varphi\|+1)+C_5\|\varphi\left(X_0\right)\|_q\right)^2,\nonumber
\end{eqnarray}
where $C_6=C\kappa^2C_4C_2^2$.

\bigskip

\textsl{Third case}: $j, k\leq i$.

\medskip

Again, using three times the operator properties, 
\begin{eqnarray}
&&|\esp\left(\varphi\left(X_0\right)\varphi\left(X_i\right)\varphi\left(X_{i+j}\right)\varphi\left(X_{i+j+k}\right)\right)|\nonumber\\
&&=|\esp\left(\varphi\left(X_0\right)\esp\left(\varphi\left(X_i\right)\esp\left(\varphi\left(X_{i+j}\right)\esp\left(\varphi\left(X_{i+j+k}\right)|\F_{i+j}\right)|\F_{i}\right)|\F_{0}\right)\right)|\nonumber\\
&&=|\esp\left[\varphi\left(X_0\right)\left[\esp\left(\varphi\left(X_i\right)\esp\left(\varphi\left(X_{i+j}\right)\esp\left(\varphi\left(X_{i+j+k}\right)|\F_{i+j}\right)|\F_{i}\right)|\F_{0}\right)\right.\right.\nonumber\\
&&\left.\left.\qquad-\esp\left(\varphi\left(X_i\right)\varphi\left(X_{i+j}\right)\varphi\left(X_{i+j+k}\right)\right)\right]\right]|\nonumber\\
&&\leq \|\varphi\left(X_0\right)\|_q\|P^i\left(\varphi P^j\left(\varphi P^k\varphi\right)\right)\left(X_0\right)-\Pi \left(\varphi P^j\left(\varphi P^k\varphi\right)\right)\left(X_0\right)\|_p\nonumber\\
&&\leq \|\varphi\left(X_0\right)\|_qC\kappa\theta^i\|\varphi P^j\left(\varphi P^k\varphi\right)\|\label{cas3.1}
\end{eqnarray}
and
\begin{eqnarray}
\|\varphi P^j\left(\varphi P^k\varphi\right)\|&\leq& \|\varphi\|\| P^j\left(\varphi P^k\varphi\right)\|\nonumber\\
&\leq&\|\varphi\|\left(\| P^j\left(\varphi P^k\varphi\right)-\Pi \varphi P^k\varphi\|+|\esp\left(\varphi\left(X_0\right)\varphi\left(X_k\right)\right)|\right)\nonumber\\
&\leq & \|\varphi\|\left(\kappa\theta^j\|\varphi P^k\varphi\|+|\esp\left(\varphi\left(X_0\right)\varphi\left(X_k\right)\right)|\right)\nonumber\\
&\leq& \kappa^2\theta^{j+k}\|\varphi\|^3+C\kappa\|\varphi\left(X_0\right)\|_q\theta^k\|\varphi\|^2,\label{cas3.2}
\end{eqnarray}
where we used inequality $\left(\ref{cas2.3}\right)$ at the last line.

\medskip

From $\left(\ref{cas3.1}\right)$ and $\left(\ref{cas3.2}\right)$, we derive
\begin{eqnarray}
&&|\esp\left(\varphi\left(X_0\right)\varphi\left(X_i\right)\varphi\left(X_{i+j}\right)\varphi\left(X_{i+j+k}\right)\right)|\nonumber\\
&&\leq C\kappa^3\|\varphi\left(X_0\right)\|_q\theta^{i+j+k}\|\varphi\|^3+C^2\kappa^2\|\varphi\left(X_0\right)\|_q^2\theta^{i+k}\|\varphi\|^2\nonumber\\
&&\leq C_7\theta^i\|\varphi\left(X_0\right)\|_q\|\varphi\|^2\left(\|\varphi\|+\|\varphi\left(X_0\right)\|_q\right)\nonumber
\end{eqnarray}
where $C_7=\max\{C\kappa^3, C^2\kappa^2\}$.

\medskip

With this last inequality and $\left(\ref{cas1.1}\right)$, the sum $\left(\ref{i}\right)$ can be bounded in the same way as before. We use the integer $n_0$ to get
\begin{eqnarray}
&&\sum_{i=1}^n\sum_{j,k\leq i}\esp\left(\varphi\left(X_0\right)\varphi\left(X_i\right)\varphi\left(X_{i+j}\right)\varphi\left(X_{i+j+k}\right)\right)\nonumber\\
&&\leq \sum_{i=1}^{3n_0-3}i^2\|\varphi\left(X_0\right)^4\|_{1}+C_7\sum_{i=3n_0-2}^n i^2\theta^i\|\varphi\left(X_0\right)\|_q\|\varphi\|^2\left(\|\varphi\|+\|\varphi\left(X_0\right)\|_q\right)\nonumber\\
&&\leq 27\left(n_0-1\right)^3\|\varphi\left(X_0\right)^4\|_{1}+ C_7\left(\|\varphi\left(X_0\right)\|_q+\|\varphi\left(X_0\right)\|_q^2\right)\sum_{i=3n_0-2}^ni^2\theta^{i-3n_0}.\nonumber
\end{eqnarray}

The sum is bounded by the corresponding series which is finite (majoration by $C_8(n_0-1)^2$, where $C_8$ depends only on $\theta$). 
So, we can conclude the study of the third case by
\begin{eqnarray}
&&\sum_{i=1}^n\sum_{j,k\leq i}\esp\left(\varphi\left(X_0\right)\varphi\left(X_i\right)\varphi\left(X_{i+j}\right)\varphi\left(X_{i+j+k}\right)\right)\nonumber\\
&&\leq 27C_2^3\|\varphi\left(X_0\right)^4\|_{1}\log^3(\|\varphi\|+1)\label{cas3}\\
&&\quad+ C_7C_8C_2^2\left(\|\varphi\left(X_0\right)\|_q+\|\varphi\left(X_0\right)\|_q^2\right)\log^2(\|\varphi\|+1).\nonumber
\end{eqnarray}

\bigskip

To conclude let $K$ be the maximum of all the constants appearing in $\left(\ref{cas1}\right)$, $\left(\ref{cas2}\right)$ and $\left(\ref{cas3}\right)$:

\begin{eqnarray*}
&&\esp\left[S_n\left(\varphi\right)^4\right]\\
&&\leq 4!K\left[n\|\varphi\left(X_0\right)^4\|_1\log^3(\|\varphi\|+1)\right.\\
&&\quad\quad\quad+n\left(\|\varphi\left(X_0\right)^3\|_q+\|\varphi\left(X_0\right)^2\|_q+\|\varphi\left(X_0\right)\|_q+\|\varphi\left(X_0\right)\|_q^2\right)\log^2(\|\varphi\|+1)\\
&&\quad\quad\quad\left.+n^2\left(\|\varphi\left(X_0\right)^2\|_1\log(\|\varphi\|+1)+\|\varphi\left(X_0\right)\|_q\right)^2\right].
\end{eqnarray*}
\cqfd

\section{Fourth moment inequality for dynamical systems}\label{s4}

In Section \ref{s2}, we dealt with homogeneous Markov chains through their operators. As usual the techniques can be applied to dynamical systems, using transfer operator.
Here we state the result for dynamical systems but the proof (which is essentially the same as in Section 3) is left to the reader.

\medskip

Let $(\Omega,\A,\mu)$ be a probability space and $T$ a measurable measure preserving transformation (i.e.
$\forall A\in\A, \mu(T^{-1}A)=\mu(A)$).
Let us consider the Perron-Frobenius operator (or the transfer operator) of $T$, $P:\L^1(\mu)\longrightarrow\L^1(\mu)$
defined by
$$
\int_{\Omega}Pf(x)g(x)\dif\mu(x)=\int_{\Omega}f(x)g\circ T(x)\dif\mu(x)
$$
for all $f\in\L^1(\mu)$ and $g\in\L^\infty(\mu)$.

As in the Markov case, we assume there exists a Banach space $(\B,\|.\|)$ of $\mu$-measurable functions from $\Omega$ to $\R$ which
contains $\bf{1}$ and satisfies (\ref{algebra}) and that $P$ verifies same assumptions:
\begin{enumerate}
\item[(i)] there exist $\kappa>0$ and $0<\theta<1$ such that for all $f\in\B$, 
$$
\|P^nf-\Pi f\|\leq \kappa\theta^n\|f\|
$$
where $\Pi f=\esp_\mu\left(f\right)\bf{1}$.

\item[(ii)] there exists $p\geq 1$ such that $\left(\B,\|. \|\right)$ is continuously included in $\left(\L^p\left(\mu\right), \|.\|_p\right)$.
\end{enumerate}

Again, it corresponds to some quasi-compactness of the Perron-Frobenius operator, see Hennion and Herv\'e \cite{HenHer01}, Baladi \cite{Bal00}.

Under these assumptions, we have same fourth moment bound: 

\begin{thm}\label{thm2}
For all $f\in\B$ such that $\esp_\mu\left(f\right)=0$, $f$ is bounded and
${\rm m}_f=\max\{1,\sup_x|f(x)|\}$,
$$
\esp_\mu\left[\left(\sum_{i=1}^nf\circ T^i\right)^4\right]\leq
K{\rm m}_f^3\left[n\|f\|_q\log^3(\|f\|+1)
+n^2\|f\|_q^2\log^2(\|f\|+1)\right].
$$
where $\frac{1}{p}+\frac{1}{q}=1$.
\end{thm}

\section{Applications}\label{s5}

We give some examples where the fourth moment inequality applies and then leads to some empirical process invariance principles.

\paragraph{Uniform ergodicity.}

Let $(X_n)_{n\ge 0}$ be a Markov chain on the state space $E$.
Denote by $(\B^\infty,\|.\|_\infty)$ the space of bounded measurable functions from $E$ to $\R$ provided with the uniform norm.
One says that the Markov chain $(X_n)_{n\ge 0}$ is {\it uniformly ergodic} if it is $\B^\infty$-geometrically ergodic.
This condition is equivalent to the fact that the process satisfies the so-called Doeblin's condition (see Meyn and Tweedie
\cite{MeyTwe93}).

In this situation, if $X_0$ has a distribution function which is enough regular,
our fourth moment bound (Corollary \ref{cor}) implies inequality (\ref{mb}). Then tightness follows and empirical process invariance principle will follow from the multivariate central limit theorem.
Note that this result is already known since Billingsley \cite{Bil68}.

See many examples of uniformly ergodic Markov chains in Meyn and Tweedie \cite{MeyTwe93}, like T-chains on compact spaces.
Another example is given by the Knudsen Gas model (see e.g. P\`ene \cite{Pen05}).

\paragraph{Expanding maps.}

For a large class of expanding maps, empirical process invariance principles have already been established in Collet, Martinez and Schmitt \cite{ColMarSch04} and Dedecker and Prieur \cite{DedPri07}. 
Expanding maps is an example of dynamical system on an interval. 
Such transformations are studied in Broise 
\cite{Bro96}, as continuous fraction expansion, $\beta$-transformation, Gauss map. In many cases, one can show that
the Perron-Frobenius operator admits a spectral gap on the space BV of bounded variation functions. Since the indicator functions belong to BV, Theorem \ref{thm2} with condition on the initial distribution prove inequality (\ref{mb}).

Gou\"ezel \cite{Gou08} gave an example of an expanding map of the interval which has a spectral gap on the space
of Lipschitz functions but not on BV. For Gou\"ezel's example, Theorem \ref{thm2} holds on the space of Lipschitz functions.

\paragraph{Subshifts.}

Let $E$ be a finite set and $\mathcal{E}=E^\N$. The metric $d$ defined on $\mathcal{E}$ is 
$$
d(x,y)=2^{-\inf\{k\ge0\,:\,x_k\neq y_k\}}.
$$
Let $A=(a(i,j))_{i,j\in E}$ be a matrix with coefficients in $\{0,1\}$ and
$$
\Omega=\{x\in\mathcal{E}\,:\,a(x_k,x_{k+1})=1,\forall k\ge 0\}.
$$
Write $T$ the shift operator on $\Omega$, i.e. $(Tx)_k=x_{k+1}$, $\forall k\ge 0$.
Denote by $\B$ the space of complex valued functions on $\Omega$, which are Lipschitz continuous with respect to
the metric $d$. The norm on $\B$ is $\|.\|=\|.\|_\infty+m(.)$, where
\begin{equation}\label{m}
m(f)=\sup\left\{\frac{|f(x)-f(y)|}{d(x,y)}, x\ne y\right\}.
\end{equation}
Note that $(\Omega, d)$ is compact and then $\B\subset\L^\infty$.
The Ruelle-Perron-Frobenius Theorem shows that the transfer operator $P$ has some quasi-compact properties on $\B$.
See Parry and Pollicott \cite{ParPol90} or Hennion and Herv\'e \cite{HenHer01}. If $1$ is the only eigenvalue of modulus $1$ and is simple, then
conditions (i) and (ii) holds and Theorem \ref{thm2} is satisfied. If $f$ is a Lipschitz continuous function on $\Omega$,
then by Dehling, Durieu and Voln\'y \cite{DehDurVol08}, an empirical process invariance principle is satisfied for the process
$(f\circ T^i)_{i\ge0}$.

\paragraph{Linear processes.}

Let $(A,\|.\|_A)$ be a separable Banach space.
Let $(a_i)_{i\ge 0}$ be a sequence of linear form on $A$ such that $\sum_{i\ge 0}|a_i|<\infty$, where
$|a|=\sup\{|a(x)|,\|x\|_A\le 1\}$.
Let $(\xi_i)_{i\in\Z}$ be an i.i.d. sequence of bounded $B$-valued random variables, where $B$ is a compact subset of $A$.
We define the $\R$-valued linear process
$$
X_k=\sum_{i\ge 0}a_i(\xi_{k-i}).
$$
If $A$ is a finite set, linear processes can be view as functionals of subshifts.
Here, in the general case, we use a slightly different metric. 
Assume there exists $\rho<1$ and $C>0$ such that $|a_i|\le C\rho^i$, for all $i\ge0$.
We defined on $B^\N$ the metric
$$
d(x,y)=\sum_{i\ge0}\rho^i\|x_i-y_i\|_A,
$$
where $x=(x_i)_{i\ge0}$ and $y=(y_i)_{i\ge0}$.
Denote by $Y_k=(\dots,\xi_{-1},\xi_0)$. Then $(Y_k)_{k\ge 0}$ is a Markov chain on $B^\N$
and
one can shows that $(Y_k)_{k\ge 0}$ satisfies (\ref{exp}) on the space $\B$ of Lipschitz continuous functions from $B^\N$ to $\R$.
Indeed, if $f\in\B$
\begin{eqnarray*}
|\esp(f(Y_k)|Y_0=y)-\esp(f(Y_k)|Y_0=x)|&=&
|\esp(f(\dots,y_0,\xi_1,\dots,\xi_k)-f(\dots,x_0,\xi_1,\dots,\xi_k))|\\
&\le& C\rho^k\|y-x\|_A
\end{eqnarray*}
and then (\ref{exp}) holds.
Hence, $(Y_k)_{k\in\Z}$ is strongly ergodic with respect to the space of Lipschitz functions on $B^\N$ and Theorem
\ref{thm} holds. It is clear that $f(x)=\sum_{i\ge0}a_i(x_i)$ is a Lipschitz function (on $B^\N$) and then fourth moment bound
holds for $(X_k)_{k\ge 0}$ on the space of Lipschitz functions (on $\R$).

\paragraph{Random iterative Lipschitz models.}

Let $G$ be a semi-group of Lipschitz transformations of a metric space $(E,d)$ and let $\mathcal{G}$ be a 
$\sigma$-algebra on $G$. We assume that the action of $G$ on $E$ is measurable.

Let $(g_n)_{n\ge 1}$ be an i.i.d. sequence of $G$-valued random variables with distribution $\eta$. Let $X_0$ be a $E$-valued random variable independent of $(g_n)_{n\ge 1}$. We consider the Markov chain $(X_n)_{n\ge 0}$ defined by
$$
X_n=g_n(X_{n-1})
$$
with transition operator $P$ definied by
$$
Pf(x)=\int_Gf(g(x))\dif\eta(g).
$$

One says $\eta$ is contracting if 
$$
\lim_{n}\sup\left\{\int_G\frac{d(g(x),g(y))}{d(x,y)}\dif \eta^{*n}(g)\,:\, x,y\in E, x\ne y\right\}^{\frac{1}{n}}<1
$$
where $\eta^{*n}$ denotes the distribution of $g_n\circ\dots\circ g_1$.

\medskip

Assume $(E,d)$ is compact. Let $\B_0$ be the space of $\mathbb{C}$-valued Lipschitz continuous functions
on $E$ provided with the norm $\|.\|_0=\|.\|_\infty+m(.)$, where $m$ is defined as in (\ref{m}).
It is shown in Hennion and Herv\'e \cite{HenHer01}, Thm.X.3, that if $\eta$ is contracting, then there exists a unique $P$-invariant measure on $E$ and $(X_n)_{n\ge 0}$ is $\B_0$-geometrically ergodic with respect to this measure.
Then fourth moment bound holds and thanks to Dehling, Durieu and Voln\'y \cite{DehDurVol08}, an empirical invariance principle follows.
One example of application is given by products of invertible random matrices (see Hennion and Herv\'e \cite{HenHer01} X.4).

\medskip

In the case where $(E,d)$ is not compact but every closed ball in $E$ is compact, one can have a similar result but 
with another Banach space (see Hennion and Herv\'e \cite{HenHer01} Thm.X.4). Here, the Banach space is the space $\B_1$ of locally 
Lipschitz functions with weight. These are the $\mathbb{C}$-valued functions $f$ such that
$$
m_1(f)=\sup\left\{\frac{|f(x)-f(y)|}{d(x,y)p(x)p(y)}\,:\,x\ne y\right\}<\infty
$$
where $p(x)=1+d(x,x_0)$ for a fixed $x_0\in E$. The norm is
$$
\|f\|_1=\sup\left\{\frac{|f(x)|}{p(x)^2}\,:\, x\in E\right\}+m_1(f).
$$
As an example, we can mention a large class of autoregressive models.

\paragraph{Autoregressive models.}

The process $(X_n)_{n\ge 0}\subset\R^d$ is called autoregressive with initial value
$X_0\in\R^d$ if it satisfies, for all $n\ge 1$,
$$
X_n=AX_{n-1}+Y_n
$$
where $A\in\mathcal{M}(\R^d)$ and $(Y_n)_{n\ge 1}\subset\R^d$ is an i.i.d. sequence of random variables, independent of $X_0$.
See Hennion and Herv\'e \cite{HenHer01} Thm.X.16 for conditions under which $(X_n)_{n\in\N}$ is $\B_1$-geometrically ergodic.

\small

\bibliographystyle{plain}

\end{document}